\documentclass{gtart_h}

\def\ifplaintex{\expandafter\ifx\csname documentclass\endcsname\relax}

\ifplaintex 
\else
\fi

\def\gt{{\mathsurround=0pt\it $\cal G\mskip-2mu$eometry \&\ 
$\cal T\!\!$opology}}        

\def\gtm{{\mathsurround=0pt\it $\cal G\mskip-2mu$eometry \&\ 
$\cal T\!\!$opology $\cal M\mskip-1mu$onographs}}    

\def\gtp{{\mathsurround=0pt\it $\cal G\mskip-2mu$eometry \&\ 
$\cal T\!\!$opology $\cal P\!$ublications}}  

\def\recd{{\small Received:\qua\receiveddate\ifx\reviseddate\relax
\else\qquad Revised:\qua\reviseddate\fi\par}} 

\def\volumenumber#1{\def\thevolumenumber{#1}}
\def\volumeyear#1{\def\thevolumeyear{#1}}
\def\volumename#1{\def\thevolumename{#1}}
\def\papernumber#1{\def\thepapernumber{#1}}
\def\pagenumbers#1#2{\def\startpage{#1}\def\finishpage{#2}}
\def\published#1{\def\publishdate{#1}}
\def\received#1{\def\receiveddate{#1}}
\def\revised#1{\def\reviseddate{#1}}
\def\accepted#1{\def\accepteddate{#1}}
\def\asciititle#1{\def\theasciititle{#1}}

\long\def\asciiabstract#1{\long\def\theasciiabstract{#1}}

\let\\\par
\let\thevolumenumber\relax\let\thepapernumber\relax
\let\thevolumeyear\relax\let\startpage\relax
\let\finishpage\relax\let\publishdate\relax\let\receiveddate\relax
\let\reviseddate\relax\let\accepteddate\relax\let\theasciititle\relax
\let\theasciiauthors\relax
\let\theasciiabstract\relax

\let\theerratum\relax\let\theasciiemail\relax
\let\theshortauthors\relax\let\theshorttitle\relax

\def\startpage{1}\def\finishpage{15}\def\thepapernumber{77}

\volumenumber{2}
\volumename{Proceedings of the Kirbyfest}
\volumeyear{1999}

\long\def\maketitlep{   

\count0=\startpage

\gtm\nl        
{\small Volume \thevolumenumber: \thevolumename\nl 
\ifx\theerratum\relax\else Erratum \erratumnumber\nl\fi
Pages \startpage--\finishpage\nl}

\vglue 0.1truein   

{\parskip=0pt\leftskip 0pt plus 1fil\def\\{\par\smallskip}{\ifplaintex\large
\else\Large\fi\bf\thetitle}\par\medskip}   
\vglue 0.05truein 

{\parskip=0pt\leftskip 0pt plus 1fil\def\\{\par}{\sc\theauthors}
\par\medskip}%
 
\vglue 0.03truein 

{\small\leftskip 25pt\rightskip 25pt{\bf Abstract}\stdspace\theabstract

{\bf AMS Classification}\stdspace\theprimaryclass
\ifx\thesecondaryclass\relax\else; \thesecondaryclass\fi\par
{\bf Keywords}\stdspace \thekeywords\par}\vglue 7pt

}   

\font\phead=cmsl9 scaled 950
\font=cmsl9 scaled 1050
\font\pnum=cmbx10 scaled 913
\font\lnum=cmbx10 
\font\pfoot=cmsl9 scaled 950
\font=cmsl9 scaled 1050
\ifplaintex
\headline{\vbox to 0pt{\vskip -4.5mm\line{\small\phead\ifnum
\count0=\startpage ISSN 1464-8997 (on line)
1464-8989 (printed) \hfill {\pnum\folio}\else\ifodd\count0\def\\{ }%
\ifx\theshorttitle\relax\thetitle\else\theshorttitle\fi\hfill{\pnum\folio}
\else\def\\{ and }{\pnum\folio}\hfill\ifx\theshortauthors\relax\theauthors
\else\theshortauthors\fi\fi\fi}\vss}}
\footline{\vbox to 0pt{\vglue 0mm\line{\small\pfoot\ifnum\count0=\startpage
Published \publishdate:\qua\copyright\ \gtp\hfill\else
\gtm, Volume \thevolumenumber\ (\thevolumeyear)\hfill\fi}\vss
}}
\else
\makeatletter
\def\@oddhead{{\small\lhead\ifnum\count0=\startpage ISSN 1464-8997 (on line)
1464-8989 (printed) \hfill {\lnum\number\count0}\else\ifodd\count0
\def\\{ }\ifx\theshorttitle\relax \thetitle \else\theshorttitle\fi\hfill
{\lnum\number\count0}\else\def\\{ and }{\lnum\number\count0}
\hfill\ifx\theshortauthors\relax 
\theauthors\else\theshortauthors\fi\fi\fi}}\def\@evenhead{@oddhead}
\def\@oddfoot{\small\lfoot\ifnum\count0=\startpage Published \publishdate:\qua\copyright\ \gtp\hfill\else
\gtm, Volume \thevolumenumber\ (\thevolumeyear)\hfill\fi}
\def\@evenfoot{@oddfoot}
\makeatother
\fi

\let\maketitlepage\maketitlep
\let\makeshorttitle\maketitlepage
\let\maketitle\maketitlepage

\newwrite\gtoutfile
\long\gdef\makeheadfile{  
{\def\\{, }\def\s{ }
\immediate\openout\gtoutfile head.xxx
\immediate\write\gtoutfile{Proxy-for: \ifx\theasciiauthors\relax
\theauthors\else\theasciiauthors\fi\s<\ifx\theasciiemail\relax\theemail\else\theasciiemail\fi>}
\immediate\write\gtoutfile{\noexpand\\}
\immediate\write\gtoutfile{Authors: \ifx\theasciiauthors\relax
\theauthors\else\theasciiauthors\fi}
{\def\\{ }\immediate\write\gtoutfile{Title: \ifx\theasciititle\relax
\thetitle\else\theasciititle\fi}}
\immediate\write\gtoutfile{Subj-class: GT or SG, GR etc}
\immediate\write\gtoutfile{MSC-class: \theprimaryclass\ifx\thesecondaryclass\relax\else, \thesecondaryclass\fi}
\immediate\write\gtoutfile{Journal-ref: Geom. Topol. Monogr. \thevolumenumber\s
(\thevolumeyear) \startpage-\finishpage}
\immediate\write\gtoutfile{Comments: Published by Geometry and Topology Monographs at}
\immediate\write\gtoutfile{\s\s\s  http://www.maths.warwick.ac.uk/gt/GTMon\thevolumenumber/paper\thepapernumber.abs.html}
\immediate\write\gtoutfile{\noexpand\\}
\immediate\write\gtoutfile{}
\ifx\theasciiabstract\relax
\immediate\write\gtoutfile{\theabstract}\else
\immediate\write\gtoutfile{\theasciiabstract}\fi
\immediate\write\gtoutfile{}
\immediate\write\gtoutfile{\noexpand\\}
\immediate\write\gtoutfile{}
\immediate\closeout\gtoutfile}}  

\def\maketitlepage{\maketitlep\makeheadfile}
\let\makeshorttitle\maketitlepage
\let\maketitle\maketitlepage

\volumenumber{7}
\volumename{Proceedings of the Casson Fest}
\volumeyear{2004}
\papernumber{16}
\pagenumbers{493}{507}
\received{9 September 2003}
\revised{25 May 2004}
\accepted{10 May 2004}
\published{13 December 2004}

\usepackage{amsmath, amssymb, amscd}

\theoremstyle{plain}
\newtheorem{thm}{Theorem}[section]
\newtheorem{lem}[thm]{Lemma}

\newtheorem{conj}[thm]{Conjecture}

\theoremstyle{definition}
\newtheorem*{defn}{Definition}

\newcommand{\Cplx}{{\mathbb{C}}}

\newcommand{\Zed}{{\mathbb{Z}}}

\newcommand{\BM}{{\mathrm{BM}}}

\newcommand{\Hecke}{{\mathcal{H}}}
\newcommand{\Local}{{\mathcal{L}}}

\begin{document}
\title{Homological representations\\of the Iwahori--Hecke algebra}
\asciititle{Homological representations\\of the Iwahori--Hecke algebra}
\author{Stephen Bigelow}
\address{Department of Mathematics,
         University of California at Santa Barbara\\
         California 93106, USA}
\email{bigelow@math.ucsb.edu}

\begin{abstract}
Representations of the Iwahori--Hecke algebra of type $A_{n-1}$
are equivalent to representations of the braid group $B_n$
for which the generators satisfy a certain quadratic relation.
We show how to construct such representations
from the natural action of $B_n$
on the homology of configuration spaces of the punctured disk.
We conjecture that all irreducible representations
of $\Hecke_n$ can be obtained in this way,
even for non-generic values of $q$.
\end{abstract}

\asciiabstract{%
Representations of the Iwahori-Hecke algebra of type A_{n-1} are
equivalent to representations of the braid group B_n for which the
generators satisfy a certain quadratic relation.  We show how to
construct such representations from the natural action of B_n on the
homology of configuration spaces of the punctured disk.  We conjecture
that all irreducible representations of Hecke_n can be obtained in
this way, even for non-generic values of q.}



\maketitle\vspace{-5pt}

\leftskip 25pt {\small\it Dedicated to Andrew Casson,
whose commitment to accuracy and elegance
has greatly influenced my work}\rightskip 25pt

\leftskip 0pt\section{Introduction}\rightskip 0pt

This paper is concerned with the representation theory of
the braid group $B_n$ and the Iwahori--Hecke algebra $\Hecke_n$.
The simplest and least motivated definition of these objects
is by generators and relations.
The braid group $B_n$ is
the group generated by $\sigma_1,\dots,\sigma_{n-1}$
with defining relations
\begin{itemize}
\item $\sigma_i \sigma_j = \sigma_j \sigma_i$ if $|i-j| > 1$,
\item $\sigma_i \sigma_j \sigma_i = \sigma_j \sigma_i \sigma_j$ if $|i-j|=1$.
\end{itemize}
Let $q$ be a unit in an integral domain $K$.
The Iwahori--Hecke algebra $\Hecke_n(K,q)$, or simply $\Hecke_n$,
is the $K$--algebra generated by $T_1,\dots,T_{n-1}$ with defining relations
\begin{itemize}
\item $T_i T_j = T_j T_i$ if $|i-j| > 1$,
\item $T_i T_j T_i = T_j T_i T_j$ if $|i-j| = 1$,
\item $(T_i + 1)(T_i - q) = 0$.
\end{itemize}

The Iwahori--Hecke algebra
plays a central role in representation theory.
If $q=1$ then $\Hecke_n$ is the group algebra
$KS_n$ of the symmetric group $S_n$.
If $K$ has characteristic zero and $q$ is {\em generic}
then $\Hecke_n$ is isomorphic to $KS_n$.
Here, generic means $q^i \neq 1$ for $i=2,\dots,n$.
Most research on $\Hecke_n$ is now concerned with
understanding the non-generic case.

There is an obvious map from $B_n$ to $\Hecke_n$
given by $\sigma_i \mapsto T_i$.
Thus the representations of $\Hecke_n$
are equivalent to those representations of $B_n$
for which the generators satisfy
$(\sigma_i + 1)(\sigma_i - q) = 0$.
The aim of this paper is to describe
a topological method of constructing such representations.

The main construction is given in Section \ref{sec:representation}.
Briefly, the idea is as follows.
A braid can be represented by
a homeomorphism from the punctured disk $D_n$ to itself.
This induces a homeomorphism from
a configuration space $C_m(D_n)$ to itself.
This in turn induces an automorphism of a certain module
obtained from $C_m(D_n)$ using homology theory.

The homology module is defined
using a local coefficient system
that depends on a representation of $B_m$.
The main result of this paper, Theorem \ref{thm:hecke},
states that if the representation of $B_m$
satisfies the required quadratic relation
then so does the resulting representation of $B_n$,
up to some rescaling.

In preparation for this result,
Section \ref{sec:HmBM} studies a slightly simpler homological representation,
and Section \ref{sec:hermitian} uses Poincar\'e duality
to define a Hermitian form that is preserved by the action of $B_n$.
Section \ref{sec:conclusion} gives a conjectural method
to obtain all irreducible representations of Iwahori--Hecke algebras
by starting from the trivial representation and inductively
applying the construction from Section \ref{sec:representation}.
We conclude with some examples in support of this conjecture,
and some more open questions motivated by it.

Lawrence \cite{rL96} has some results
very similar to those in this paper.
The construction and the outcome
are almost identical to mine,
although the method of proof is quite different.
There are two main advantages to the approach in this paper.
First,
it explicitly identifies the representation
as an image of a natural map from homology to relative homology,
whereas \cite{rL96} uses a less well motivated quotient.
Second,
it works even when $q$ is a root of unity,
whereas \cite{rL96} needs to assume $q$ is generic.
On the other hand,
\cite{rL96} has the advantage of
identifying precisely what representations are obtained,
whereas this paper only does so conjecturally
in Conjecture \ref{conj:d}.

We use the following notation throughout this paper.
\begin{itemize}
\item $D$ is the unit disk centered at $0$ in the complex plane.
\item $-1 < p_1 < \dots < p_n < 1$ are distinct points on the real line.
\item $D_n$ is the $n$--times punctured disk $D \setminus \{p_1,\dots,p_n\}$.
\end{itemize}
Then $B_n$ is the group of homeomorphisms $f \co D_n \to D_n$
such that $f | \partial D$ is the identity,
taken up to isotopy relative to $\partial D$.
The generator $\sigma_i$ corresponds to
a homeomorphism that exchanges $p_i$ and $p_{i+1}$
by a counterclockwise half twist.

We also use the definition of $B_n$
as the fundamental group of a configuration space.
If $X$ is a surface,
let $C_n(X)$ denote the configuration space
of unordered $n$--tuples of distinct points in $X$.
Then $B_n = \pi_1(C_n(D),\{p_1,\dots,p_n\})$.
The generator $\sigma_i$ corresponds to a loop
in which $p_i$ and $p_{i+1}$
switch places by a counterclockwise half twist,
while the other points remain fixed.

The main results of this paper require $K$ to be
a field with a conjugation operation.
However the construction in Section \ref{sec:representation}
and the results in Section \ref{sec:HmBM}
apply when $K$ is any integral domain.

\rk{Acknowledgement}This research was partly supported by NSF grant
DMS-0307235 and Sloan Fellowship BR-4124.

\section{Constructing a representation}
\label{sec:representation}

In this section we describe a construction that
produces a representation $W$ of $B_n$
from a representation $V$ of $B_m$,
where $n$ and $m$ are any non-negative integers.
We leave open the possibility that $W$ is the zero vector space,
although this should not technically be called a representation of
$B_n$.
In Section \ref{sec:conclusion}
we give some examples and a conjecture
suggesting that in interesting cases
$W$ will be in some sense more sophisticated than $V$.

Let $B_m$ act by automorphisms on a
non-zero
$K$--module $V$.
Let $C$ be the configuration space $C_m(D_n)$.
Let $c_0 \in C$ be a set of $m$ points on $\partial D$.
We will use the action of $B_m$ on $V$
to define an action of $\pi_1(C,c_0)$ on $V$.
The construction does not require a deep understanding of $\pi_1(C,c_0)$.
However we mention that $\pi_1(C,c_0)$ can be described succinctly
as the subgroup of $B_{n+m}$ consisting of geometric braids
for which the first $n$ strands are straight lines.

Let $f \co C \to \Cplx \setminus \{0\}$ be the map
  $$f(\{z_1,\dots,z_m\}) = \prod_{i=1}^m \prod_{j=1}^n (z_i-p_j).$$
Let $w \co \pi_1(C,c_0) \to \Zed$
be the map that takes a loop $\gamma$ in $C$
to the winding number of $f \circ \gamma$ around $0$.

Let $i \co C \to C_m(D)$
be the map induced by the inclusion $D_n \subset D$.
Then $i_*$ is a surjective homomorphism from $\pi_1(C,c_0)$ to $B_m$.
Let $\pi_1(C)$ act on $V$ by
  $$g(v) = q^{w(g)} i_*(g)(v)$$
for all $g \in \pi_1(C)$ and $v \in V$.
Let $\Local$ be the $V$--bundle over $C$
whose monodromy is given by this action.

Let $H_i(C;\Local)$ be the homology module
with local coefficients given by $\Local$.
We recall briefly how this is defined.
A singular chain is a formal sum of terms of the form $n \sigma$,
where $\sigma \co \Delta^i \to C$ is a singular $i$--simplex in $C$
and $n \co \Delta^i \to \Local$ is a lifting of $\sigma$.
The obvious boundary map gives rise to a chain complex,
whose homology modules are $H_i(C;\Local)$.
We will be particularly interested in
the middle homology $H_m(C;\Local)$.

Relative homology with local coefficients
is defined in the usual way.
We also use {\em Borel--Moore} homology, defined by
  $$H_m^\BM(C;\Local) = \lim_{\leftarrow} H_m(C,C \setminus A;\Local),$$
where the inverse limit is taken over all compact subsets $A$ of $C$.

There is an action of the braid group $B_n$ on $H_m(C;\Local)$
defined as follows.
Let $f \co D_n \to D_n$ be a homeomorphism
that acts as the identity on $\partial D$.
Then $f$ induces a homeomorphism $f' \co C \to C$
that fixes $c_0$.
The induced automorphism $f'_*$ of $\pi_1(C,c_0)$
is such that $f'_*(g)(v) = g(v)$ for all $g \in \pi_1(C,c_0)$ and $v \in V$.
Thus there is a unique lift $\tilde{f}' \co \Local \to \Local$
that acts as the identity on the fiber over $c_0$.
This map induces an automorphism of $H_m(C;\Local)$.
This automorphism depends only on the isotopy class of $f$
relative to $\partial D$,
so gives a well defined action of $B_n$ on $H_m(C;\Local)$.

Let $W$ be the image of the map
from $H_m(C;\Local)$ to $H_m^\BM(C,\partial C;\Local)$
induced by inclusion.
Then $B_n$ acts on $W$.
This completes our construction of a representation of $B_n$
from a representation of $B_m$.

{\em
For improved readability,
we suppress mention of $\Local$ from now on,
writing for example $H_m(C)$ instead of $H_m(C;\Local)$.
}

\section{The structure of $H_m^\BM(C)$}
\label{sec:HmBM}

Let $q$ be a unit in an integral domain $K$,
let $V$ be a $K$--module with an action of $B_m$,
and let $C$ and $\Local$ be as in Section \ref{sec:representation}.
In this section we compute $H_m^\BM(C)$.
Our main interest will ultimately be in the representation $W$,
but $H_m^\BM(C)$ is easier to understand.

\begin{defn}
Let $\Pi$ be the set of sequences
$\pi = (\pi_1,\dots,\pi_{n-1})$ such that
$\pi_i$ is a non-negative integer and $\sum \pi_i = m$.
\end{defn}

For $\pi \in \Pi$,
let $U_\pi$ be the set of all $\{x_1,\dots,x_m\} \in C$
such that $x_1,\dots,x_m \in (p_1,p_n)$ and
  $$\sharp (\{x_1,\dots,x_m\} \cap (p_i,p_{i+1})) = \pi_i$$
for $i=1,\dots,n-1$.
This is an open $m$--ball,

\begin{lem}
\label{lem:HmBM}
$H_m^\BM(C)$ is the direct sum of $\binom{n+m-1}{m}$ copies of $V$,
namely the images of $H_m^\BM(U_\pi)$ for $\pi \in \Pi$.
\end{lem}

\begin{proof}
Let $C_0$ be the set of all $\{x_1,\dots,x_m\} \in C$
such that $x_1,\dots,x_m \in (p_1,p_n)$.
Then $C_0$ is the disjoint union of $U_\pi$ for $\pi \in \Pi$.
Thus $H_m^\BM(C_0)$ is the direct sum of $V_\pi$
over all $\pi \in \Pi$.
It remains to prove that the map
$H_m^\BM(C_0) \to H_m^\BM(C)$ induced by inclusion is an isomorphism.

\begin{defn}
For $\epsilon > 0$,
let $A_\epsilon$ be the set of all $\{z_1,\dots,z_m\} \in C$
such that $|z_i - z_j| \ge \epsilon$
for all distinct $i,j = 1,\dots,m$,
and $|z_i - p_j| \ge \epsilon$
for all $i = 1,\dots,m$ and $j = 1,\dots,n$.
\end{defn}

Note that each $A_\epsilon$ is compact,
and every compact subset of $C$ is contained in some $A_\epsilon$.
Thus it suffices to show that for all sufficiently small $\epsilon > 0$,
the map
  $$H_m(C_0,C_0 \setminus A_\epsilon) \to H_m(C,C \setminus A_\epsilon)$$
induced by inclusion is an isomorphism.

Let $D'_n \subset D_n$ be
a closed $\epsilon/2$ neighborhood of the interval $[p_1,p_n]$.
Let $C' = C_m(D'_n)$.
The map
  $$H_m(C',C' \setminus A_\epsilon) \to H_m(C,C \setminus A_\epsilon)$$
induced by inclusion is an isomorphism.
To see this,
note that the obvious homotopy shrinking $C$ to $C'$
is a homotopy through maps
from $(C,C \setminus A_\epsilon)$ to itself.

Let $U$ be the set of $\{z_1,\dots,z_m\} \in C'$
such that the real parts of $z_1,\dots,z_m$
are all distinct
and lie in $[p_1,p_n] \setminus \{p_1,\dots,p_n\}$.
Then $U$ is open and contains $A_\epsilon \cap C'$.
By excision, the map
  $$H_m(U,U \setminus A_\epsilon) \to H_m(C',C' \setminus A_\epsilon)$$
induced by inclusion is an isomorphism.

There is an obvious deformation retraction from $U$ to $C_0$
taking $\{z_1,\dots,z_m\}$ to $\{x_1,\dots,x_m\}$
where $x_i$ is the real part of $z_i$.
This is a homotopy through maps of pairs of spaces
from $(U,U \setminus  A_\epsilon)$ to itself.
Thus the map
  $$H_m(C_0,C_0 \setminus A_\epsilon) \to H_m(U,U \setminus A_\epsilon))$$
induced by inclusion is an isomorphism.
\end{proof}

\subsection{The Krammer representation}

As an example,
we now describe a homological definition
of the representation 
that
Krammer defined in \cite{dK00}
and proved to be faithful in \cite{dK02}.

Let $V = K$, let $t$ be a unit of $K$,
and let the generator of $B_2$ act on $V$ by multiplication by $-t$.
Then $C$ is the space of unordered pairs of distinct points in $D_n$
and $\Local$ is a $K$--bundle with monodromy as follows.
A loop in $C$ in which the two points switch places
by a small counterclockwise twist has monodromy $-t$.
A loop in which
one point makes a small counterclockwise loop around
one
puncture
while the other point remains fixed has monodromy $q$.

By Lemma \ref{lem:HmBM},
$H_2^\BM(C)$ is a free $K$--module of rank $\binom{n}{2}$.
For $1 \le i \le j \le n-1$,
let $U_{i,j}$ be the set of $\{x,y\} \in C$
such that $x \in (p_i,p_{i+1})$ and $y \in (p_j,p_{j+1})$.
Let $u_{i,j}$ be the element of $H_2^\BM(C)$
represented by some non-zero lift of $U_{i,j}$ to $\Local$.
These form the basis of $H_2^\BM(C)$
given by the proof of Lemma \ref{lem:HmBM}.

An alternative basis is as follows.
For $1 \le i < j \le n$,
let $T_{i,j}$ be an edge from $p_i$ to $p_j$
whose interior lies in the upper half plane.
Let $V_{i,j}$ be the set of pairs of distinct points in $T_{i,j}$.
Let $v_{i,j}$ be the element of $H_2^\BM(C)$
represented by some non-zero lift of $V_{i,j}$ to $\Local$.
With appropriate choices of local coordinates,
$$v_{i,j} = \sum_{i \le k \le l < j}u_{k,l}.$$
These form a basis for $H_2^\BM(C)$.

It is not hard to compute
the matrices for the action of $B_n$
with respect to this basis.
They are the same as
the matrices used by Krammer in \cite{dK00}
and \cite{dK02} to produce a faithful representation of $B_n$.
Thus we have a homological construction
of the ``Krammer representation'' over any ring.
This proves \cite[Conjecture 5.7]{sB03}.

\section{A Hermitian form}
\label{sec:hermitian}

For the rest of this paper, we assume the following.
\begin{itemize}
\item $K$ is a field,
\item $x \mapsto \bar{x}$ is an automorphism of $K$
      such that $\bar{\bar{x}} = x$,
\item $q \in K$ satisfies $q \bar{q} = 1$ and $q \neq 1$,
\item $\langle \cdot,\cdot \rangle_V \co V \times V \to K$
      is a non-singular Hermitian form
      that is invariant under the action of $B_m$.
\end{itemize}
An interesting example is $K = \Cplx$, $q$ is a root of unity,
and $V$ is a unitary representation of $B_m$.
We will not explicitly use the assumption $q \neq 1$,
but our constructions produce trivial results if $q = 1$.

Under these assumptions,
we will define a Hermitian form on $W$
that is invariant under the action of $B_n$.
This form is something of a relic
from my previous less efficient proofs of the results in this paper.
It now is only used to prove Lemma \ref{lem:HmCpartialC} below.
With more work it could probably be completely eliminated,
thus allowing us to work over a more general integral domain $K$.
However I still expect this form to play an important role
in a proof of Conjecture \ref{conj:d}.

Suppose $\pi  = (\pi_1,\dots,\pi_{n-1}) \in \Pi$,
that is, $\pi_i$ is a non-negative integer and $\sum \pi_i = m$.
Let
  $$f_1,\dots,f_m \co (I,\partial I) \to (D_n,\partial D)$$
be disjoint vertical edges, oriented upwards,
numbered from left to right,
so that for all $i=1,\dots,n-1$
there are $\pi_i$ edges passing between $p_i$ and $p_{i+1}$.
The map
  $$F = f_1 \times \dots \times  f_m$$
induces an embedding of the closed $m$--ball into $C$.
The set of lifts of $F$ to $\Local$
defines a subspace $V_\pi$ of $H_m(C,\partial C)$.

\begin{lem}
\label{lem:HmCpartialC}
$H_m(C,\partial C)$ is the direct sum of
$V_\pi$ for $\pi \in \Pi$.
\end{lem}

This follows immediately from Lemma \ref{lem:HmBM}
and Poincar\'e duality.
We therefore devote the rest of this section
to reviewing standard duality results from homology theory,
since these may be less familiar
in the context of homology with local coefficients.

Let $C$ and $\Local$ be as in Section \ref{sec:representation}.
Recall that singular cohomology $H^m(C)$
is the homology of the space of singular cochains.
A singular cochain is a function $\phi$ that assigns
to each singular cochain $\sigma \co \Delta^i \to C$
a lift $\phi(\sigma) \co \Delta^i \to \Local$.

There are versions of Poincar\'e duality
and the universal coefficient theorem using local coefficients.
Poincar\'e duality is proved
using a cap product with a fundamental class
of the usual singular homology $H_{2m}(C;K)$.
The universal coefficient theorem gives an anti-isomorphism
from cohomology to the dual of homology.
It is proved using the anti-isomorphism from $V$ to its dual
given by the effective Hermitian form on $V$.

These results give rise to forms
  $$\langle\cdot,\cdot\rangle_1 \co H_m^\BM(C) \times H_m(C,\partial C) \to K,$$
  $$\langle\cdot,\cdot\rangle_2 \co H_m(C) \times H_m^\BM(C,\partial C) \to K,$$
that are non-singular, sesquilinear, and invariant under the action of $B_n$.

Define a form
  $$\langle \cdot,\cdot \rangle_0 \co H_m(C) \times H_m(C) \to K$$
as follows.
Suppose $u,v \in H_m(C)$.
Let $u_1 \in H_m^\BM(C)$, $v_1 \in H_m(C,\partial C)$,
and $v_2 \in H_m^\BM(C,\partial C)$
be the images of $u$ and $v$ under maps induced by inclusion.
Let
  $$\langle u,v \rangle_0 = \langle u_1,v_1 \rangle_1.$$
This is equivalent to
  $$\langle u,v \rangle_0 = \langle u,v_2 \rangle_2.$$
This is Hermitian, but not necessarily non-singular.

In practice, $u,v \in H_m(C)$ are often represented by
immersed closed orientable $m$--manifolds in $C$,
together with lifts $\tilde{u}$, $\tilde{v}$ to $\Local$.
In this case $\langle u,v \rangle$ can be computed as follows.
Let the $m$--manifolds be intersect transversely in $C$.
For each point $x$ of intersection,
let $u_x$ and $v_x$ be the corresponding points
of $\tilde{u}$ and $\tilde{v}$ in the fiber over $x$.
Then $\langle u,v \rangle_0$ is the sum of
$\langle u_x,v_x \rangle_{\Local_x}$ over all points $x$ of intersection.

Recall that $W$ is the image of the map
from $H_m(C)$ to $H_m^\BM(C,\partial C)$.
Define
  $$\langle \cdot,\cdot \rangle_W \co W \times W \to K$$
as follows.
Suppose $w_1,w_2 \in W$ are the images of $u_1,u_2 \in H_m(C)$
under the map induced by inclusion.
Then let
  $$\langle w_1,w_2 \rangle_W = \langle u_1,u_2 \rangle_0.$$
Using the properties of $\langle \cdot,\cdot \rangle_0$
that we have mentioned,
it is not hard to show that $\langle \cdot,\cdot \rangle_W$
is a well defined non-singular Hermitian form
that is invariant under the action of $B_n$.

For example, suppose $V=K=\Cplx$
and the generator of $B_2$ acts as multiplication by
$-t$,
a complex number with unit norm. 
Then $H_2^\BM(C)$ is the Krammer representation
described in Section \ref{sec:representation}.
If $q$ and $t$ are algebraically independent unit complex numbers
then it can be shown to map isomorphically onto $W$.
Thus we obtain a non-singular Hermitian form
preserved by the action of $B_n$.
In \cite{rB??},
Budney proves that this form is negative definite,
and hence the representation is unitary.

\section{A representation of the Iwahori--Hecke algebra}

Let $V$ be a finite dimensional $K$--vector space with
an action of $B_m$ that preserves an effective Hermitian form.
Let $W$ be the representation of $B_n$
constructed as in Section \ref{sec:representation}.
To avoid confusion,
denote the generators of $B_m$ by $\tau_1,\dots,\tau_{m-1}$,
and the generators of $B_n$ by $\sigma_1,\dots,\sigma_{n-1}$.

\begin{thm}
\label{thm:hecke}
If the action of $B_m$ on $V$ satisfies
  $$(\tau_i - 1)(q \tau_i + 1) = 0$$
for $i = 1,\dots,m-1$,
then the resulting action of $B_n$ on $W$ satisfies
  $$(\sigma_i - 1)(\sigma_i + q) = 0$$
for $i = 1,\dots,n-1$.
\end{thm}

\begin{proof}
It suffices to show that the action of $B_n$ on $W$
satisfies $(\sigma_1 - 1)(\sigma_1 + q) = 0$,
since the generators $\sigma_i$ are all conjugate to $\sigma_1$.
In fact, we will show that this relation holds for
the action of $B_n$ on the image of the map
  $$H_m(C,\partial C) \to H_m^\BM(C,\partial C)$$
induced by inclusion.

Suppose $\pi  = (\pi_1,\dots,\pi_{n-1}) \in \Pi$,
that is, $\pi_i$ is a non-negative integer and $\sum \pi_i = m$.
Let
  $$f_1,\dots,f_m \co (I,\partial I) \to (D_n,\partial D)$$
be disjoint vertical edges, oriented upwards,
numbered from left to right,
so that for all $i=1,\dots,n-1$
there are $\pi_i$ edges passing between $p_i$ and $p_{i+1}$.
The map
  $$F = f_1 \times \dots \times  f_m$$
induces an embedding of the closed $m$--ball into $C$.
Let $v$ be represented by a lift
$\tilde{F}$ of $F$ to $\Local$.
By Lemma \ref{lem:HmCpartialC},
such elements $v$ generate $H_m(C,\partial C)$.
We show that either
$(\sigma_1 - 1)(\sigma_1 + q) v = 0$
or the image of $v$ in $H_m^\BM(C,\partial C)$ is zero.

First consider the case $\pi_1 = 0$.
Then $\sigma_1$ is the identity on the image of $F$,
so $\sigma_1$ acts as the identity on $V_\pi$.
In particular, $(\sigma_1 - 1)v = 0$.

Next consider the case $\pi_1 = 1$.
Then $\sigma_1$ is the identity on the images of $f_2,\dots,f_n$.
Thus
  $$\sigma_1 F = \sigma_1 f_1 \times f_2 \times \dots \times f_m.$$
We can homotope the edge $\sigma_1 f_1$,
keeping it to the left of $f_2$ and keeping the endpoints on $\partial D$,
to obtain a composition of paths $g \bar{f}_1$
where $g$ goes from $f_1(0)$ to $f_1(1)$,
passing upwards between $p_2$ and $p_3$,
and $\bar{f}_1$ is $f_1$ with the reverse orientation.
Let
\begin{eqnarray*}
G &=& g \times f_2 \times \dots \times f_m, \\
F'&=& \bar{f_1} \times f_2 \times \dots \times f_m,
\end{eqnarray*}
Then $\sigma_1 v = u+v'$
where $u$ and $v'$ are represented by
appropriate lifts $\tilde{G}$ and $\tilde{F}'$ of $G$ and $F'$.

Now $F'$ is the same closed ball as $F$
but with the opposite orientation.
It remains to compare the lifts $\tilde{F}'$ and $\tilde{F}$.
Consider the loop in $\alpha$ in $C$ given by
  $$\alpha(t) = \{(g \bar{f}_1)(t), f_2(0),\dots,f_m(0)\}.$$
Let $\tilde{\alpha}$ be the lift of $\alpha$
such that $\tilde{\alpha}(0) \in \tilde{F}$.
Then $\tilde{\alpha}$
goes through $\tilde{G}$ to a point
in the intersection of $\tilde{G}$ and $\tilde{F}'$,
then through $\tilde{F}'$ to $\tilde{\alpha}(1)$.
Thus $\tilde{\alpha}(1)$
is the point in $\tilde{F}'$
corresponding to the point $\tilde{\alpha}(0)$ in $\tilde{F}$.
Now $\alpha$ has monodromy $q$, so we we conclude that $v' = -q v$.
Thus $(\sigma_1 + q)v = u$.
But $(\sigma_1 - 1)u = 0$ by the previous case.
Thus $(\sigma_1-1)(\sigma_1+q)v = 0$.

Finally, consider the case $\pi_1 \ge 2$.
Let $\epsilon > 0$
and let $A_\epsilon$ be as defined in Section \ref{sec:HmBM}.
Let $U = C \setminus A_\epsilon$.
This is the set of $\{z_1,\dots,z_m\} \in C$
such that some distinct $z_i$ and $z_j$ are
within distance $\epsilon$ of each other,
or some $z_i$ is within $\epsilon$ of a puncture.
We show that the image of $u$ in $H_m(C,U \cup \partial C)$ is zero,
and hence that the image of $u$ in $H_m^\BM(C,\partial C)$ is zero.

We can homotope $f_1$,
keeping it to the left of $f_2$ and keeping the endpoints on $\partial D$,
to obtain a composition $f f_\epsilon \bar{f}$,
where $f$ is the straight line from $-1$ to $p_1 - \epsilon/2$
and $f_\epsilon$ is a circular loop of radius $\epsilon/2$ around $p_1$.
The image of $f_\epsilon \times f_2 \times \dots \times f_m$ lies in $U$.
Thus the image of $v$ in $H_m(C,U \cup \partial C)$ is $(1-q)v'$, where
$v'$ is represented by some lift of $f \times f_2 \times \dots \times f_m$.

Next homotope $f_2$,
keeping it disjoint from $f$ and $f_3$,
and keeping the endpoints on $\partial D$,
to obtain a composition $g_1 g_\epsilon g_2$,
where $g_1$ is a straight line going right
from $\partial D$ to $p_1 - (\epsilon/2)i$,
$g_\epsilon$ is a semicircle of radius $\epsilon/2$ centered at $p_1$,
and $g_2$ is a straight line
going left from $p_1+(\epsilon/2)i$ to $\partial D$.
Let
  $$G_1 = f \times g_1 \times f_3 \times \dots \times f_m,$$
  $$G_2 = f \times g_2 \times f_3 \times \dots \times f_m.$$
Then $v' = v_1 + v_2$
where $v_i$ is represented by an appropriate lift of $G_i$.

We can simultaneously homotope $g_2$ to $\bar{f}$, and $f$ to $g_1$,
by ``pushing downwards'',
keeping one endpoint of each edge on $\partial D$
and the other in an $\epsilon$--neighborhood of $p_1$,
and keeping the two edges disjoint.
Thus $v_2$ is represented by an appropriate lift of
  $$G'_2 = g_1 \times \bar{f} \times f_3 \times \dots \times f_m.$$
This has the same image as $G_1$.
It also has the same orientation,
since the first two coordinates have switched
and the second coordinate has reversed orientation.
It remains to compare the lifts of $G_1$ and $G'_2$
that represent $v_1$ and $v_2$.

Let $\alpha$ be the path in $C$ given by
  $$\alpha(t) = \{ f(1), g_\epsilon(t), f_3(0),\dots,f_m(0) \}.$$
Let $\tilde{\alpha}$ be the lift of $\alpha$ to $\Local$
such that $\tilde{\alpha}(0) \in \tilde{G}_1$.
Then $\tilde{\alpha}(1) \in \tilde{G}_2$.
Let $\beta$ be the path in $C$ given by
  $$\beta(t) = \{ \beta_1(t), \beta_2(t), f_3(0),\dots,f_m(0) \}$$
where $\beta_1$ is a path from $h_\epsilon(1)$ down to $g(1)$
and $\beta_2$ is a path from $g(1)$ down to $h_\epsilon(0)$.
Let $\tilde{\beta}$ be the lift of $\beta$
such that $\tilde{\beta}(0) \in \tilde{G}_2$.
Then $\tilde{\beta}$ follows our homotopy
from $\tilde{G}_2$ to $\tilde{G}'_2$,
so $\tilde{\beta}(1) \in \tilde{G}'_2$.
Now $\tilde{\alpha} \tilde{\beta}$
is a lift of $\alpha \beta$
that starts in $\tilde{G}_1$
and ends at the corresponding point in $\tilde{G}'_2$.
But $\alpha \beta$ is a loop with monodromy $q \tau$.
Thus $v_2 = q \tau_1 v_1$, so
  $$v' = (1+q\tau_1)v_1.$$

We can assume $|f(t) - g_1(t)| < \epsilon$ for all $t \in I$.
Let
  $$T_1 = \{(t_1,t_2,\dots,t_m) \in I^m : t_1 > t_2\},$$
  $$T_2 = \{(t_1,t_2,\dots,t_m) \in I^m : t_1 < t_2\}.$$
Then $\tilde{G_1}|T_1$ and $\tilde{G_1}|T_2$ represent
elements $w_1$ and $w_2$ of $H_m(C,U \cup \partial C)$.
We have $v_1 = w_1 + w_2$.

Let
$$H_1 = f \times f \times f_3 \times \dots \times f_m | T_1 \to C,$$
$$H_2 = f \times f \times f_3 \times \dots \times f_m | T_2 \to C.$$
Then $w_i$ is represented by an appropriate lift $\tilde{H}_i$ of $H_i$.
The image of $H_2$ is the same as that of $H_1$,
but the orientation is reversed
since the first two coordinates have switched places.
A path from a point in $\tilde{H}_1$
to the corresponding point in $\tilde{H}_2$
is given by a lift of a loop in $C$
in which the leftmost two points switch places by a half twist.
The monodromy of this loop is $\tau_1$.
Thus $w_2 = -\tau_1 w_1$.

We conclude that the image of $v$ in $H_m^\BM(C,\partial C)$
is $(1-q)(1+q\tau_1)(1-\tau_1)w_1$.
But by assumption, the action of $B_m$ on $V$ satisfies
$(1+q\tau_1)(1-\tau_1) = 0$.
Thus the image of $v$ in $H_m^\BM(C,\partial C)$ is zero.
\end{proof}

\section{Conclusion}
\label{sec:conclusion}

We can use Theorem \ref{thm:hecke} to construct
representations of the Iwahori--Hecke algebra as follows.
Suppose $V$ is a representation of $\Hecke_m$.
Let $B_m$ act on $V$ by
  $$\tau_i \mapsto -T_i^{-1}.$$
Let $W$ be the representation of $B_n$
constructed in Section \ref{sec:representation}.
By Theorem \ref{thm:hecke},
the action of $B_n$ on $W$ factors through the map
  $$\sigma_i \mapsto q T_i^{-1},$$
and so gives a representation of $\Hecke_n$.

In their groundbreaking paper \cite{DJ86},
Dipper and James classify the irreducible representations of $\Hecke_n$.
They define a representation $D_\lambda$ of $\Hecke_n$
for every partition $\lambda$ of $n$.
They show that the non-zero $D_\lambda$
enumerate the irreducible representations of $\Hecke_n$.
Furthermore they give a simple classification
of which partitions $\lambda$ give rise to non-zero $D_\lambda$.

Let $\lambda$ be a partition of $n$.
That is, $\lambda = (\lambda_1,\dots,\lambda_k)$
where $\lambda_1 \ge \dots \ge \lambda_k > 0$
and $\lambda_1 + \dots + \lambda_k = n$.
Let $\mu$ be the partition $(\lambda_2,\dots,\lambda_k)$
of $m = n-\lambda_1$.

\begin{conj}
\label{conj:d}
Suppose $q \neq 1$.
If $V$ is the representation $D_\mu$ of $\Hecke_m$
and $W$ is the representation of $\Hecke_n$
constructed by Theorem \ref{thm:hecke}
then $W = D_\lambda$.
\end{conj}

This would make it possible to inductively construct
all irreducible representations of Iwahori--Hecke algebras from
$D_{\emptyset}$, the trivial representation of
the algebra $\Hecke_0 = K$.
The conjecture is true when $K = \Cplx$ and
$q$ is a generic unit complex number,
because in this case my construction is the same
as that of Lawrence \cite{rL96}.
We now give some examples that can be computed directly.
These confirm the conjecture in some simple non-generic cases,
and give an idea of the sort of behavior to expect in general.

\subsection{The trivial representation}

Let $\lambda$ be the partition $(n)$.
Then $\mu$ is the empty partition and $D_\mu = K$.
The configuration space $C$ consists of a single point
(namely the empty set).
Thus $W = H_0(C) = K$ with the trivial action of $B_n$.
We obtain the representation $T_i \mapsto q$ of $\Hecke_n$,
which is indeed $D_\lambda$.

\subsection{The Burau representation}

Let $\lambda$ be the partition $(n-1,1)$.
Then $\mu = (1)$, $m=1$, and $V = K$.
The configuration space $C$ is simply $D_n$.
The local coefficients $\Local$ are such that
a counterclockwise loop around a puncture has monodromy $q$.
The homology $H_1(D_n)$ is an $(n-1)$--dimensional vector space.
The induced action of $B_n$ is
the well-known homological construction of the Burau representation.

The map from $H_1(D_n)$ to $H_1^\BM(D_n,\partial D)$
is an isomorphism
except if $q$ is an $n$th root of unity.
In this case, the map has a one-dimensional kernel.
The kernel is generated by a non-zero lift of $\partial D$ to $\Local$.
Such a lift exists because the monodromy around $\partial D$
is trivial when $q^n = 1$.

Thus $W$ is either the Burau representation
or an $(n-2)$--dimensional quotient of the Burau representation.
This is indeed isomorphic to $D_\lambda$.

\subsection{The $(n-2,2)$ representation}

Suppose $n \ge 4$ and $\lambda$ is the partition $(n-2,2)$.
Then $\mu = (2)$,
and $D_\mu$ is the one-dimensional representation of $\Hecke_2$
given by $T_1 \mapsto q$.
Thus $V$ is the one-dimensional representation of $B_2$
given by $\tau_1 \mapsto -q^{-1}$.

We obtain a $K$--bundle $\Local$ over
the space $C$ of unordered pairs of distinct points in $D_n$.
A loop in $C$ in which one point goes
counterclockwise around a puncture has monodromy $q$.
A loop in which the points switch places
by a counterclockwise twist has monodromy $-q^{-1}$.

Let $W_0$ be the image of the map
from $H_2(C)$ to $H_2^\BM(C)$.
The main result of \cite{sB03}
is that $W_0$ has dimension $n(n-3)/2$
and is the ``Specht module'' corresponding to $\lambda$.
However I used the assumption that
$q$ is not a square or cube root of $1$,
which I now think I can do without.

If $q^{n-2} = 1$ then
the kernel of the map from $W_0$ to $W$ has dimension $n-1$.
A basis is given by lifts of annuli of the form
$\partial D \times (p_i,p_{i+1})$.
Such lifts exist because the loop in which one point
goes around the other point and all punctures
has monodromy $q^{n-2} = 1$.
The action of $B_n$ on this kernel
is the Burau representation.

If $q^{n-1} = 1$ then
the kernel of the map from $W_0$ to $W$
has dimension $1$.
It is the space of lifts of the surface
of pairs of distinct points in $\partial D$.

These results agree with what is known about $D_\lambda$.

\subsection{Open questions}

Assuming Conjecture \ref{conj:d} is true,
it is natural to ask the following questions.

\begin{itemize}
\item Does Conjecture \ref{conj:d} have any applications to
      the representation theory of Iwahori--Hecke algebras?
\item What is the dimension of $D_\lambda$?
\item Does Conjecture \ref{conj:d} hold when $K$ is a local ring?
\item Is there a similar construction
      for representations of the Birman--Mura\-kami--Wenzl algebra?
\end{itemize}

The first question is deliberately vague.
The representation theory of Iwahori--Hecke algebras
is a well trammeled area of study,
and it might be too optimistic to think
the topological approach in this paper
will help resolve any of the outstanding open problems.
At least I hope it will give an amusing new perspective.

The second question is a specific example
of a major open problem in representation theory.
The answer is known in the sense that
there is an unfeasibly slow algorithm to compute
the dimension of $D_\lambda$.
A great deal of more practical information is also known -
enough to make it clear that it is a very deep problem.

Viewing this problem in the light of Conjecture \ref{conj:d},
we are led to investigate
the nature of the kernel of the map
from $H_m(C)$ to $H_m^\BM(C,\partial C)$.
This kernel comes from the previous term
in a long exact sequence of relative homology.
Thus we are led to ask what is the dimension
of $H_{m+1}(\partial \bar{C})$,
where $\bar{C}$ is the set $A_\epsilon$ defined in Section \ref{sec:HmBM},
or perhaps some more sophisticated compactification of $C$.

The third question has applications to
the representation theory of the symmetric group,
which is still not well understood over non-zero characteristic.
The idea is to
first understand the representations of $\Hecke_n$
when $K$ is a ring localized at $(q-1)$,
and then understand the effect of taking the quotient by $(q-1)$.
The main difficulty if $K$ is a ring
seems to be in defining the non-singular Hermitian form.
The problem is that the universal coefficient theorem
is more complicated in the presence of torsion.

The last question is inspired by
the result due to Zinno \cite{mZ01},
which states that the Lawrence--Krammer representation
is a representation of the Birman--Murakami--Wenzl (BMW) algebra.
The BMW algebra, like the Iwahori--Hecke algebra,
is a quotient of the braid group algebra
by some relations in the generators.
By analogy with the proof of Theorem \ref{thm:hecke}
it should be possible to show that,
up to some mild rescaling,
if $V$ satisfies the relations of the BMW algebra then so does $W$.


\Addresses

\end{document}